\documentclass[english]{amsart}

\usepackage[latin1]{inputenc}

\usepackage{graphicx,rotating}

\usepackage{amsthm,amsmath,amssymb}
\theoremstyle{plain}
\newtheorem{thm}{Theorem}[section]
\newtheorem{prop}[thm]{Proposition}
\newtheorem{lem}[thm]{Lemma}
\newtheorem{conj}[thm]{Conjecture}
\newtheorem{cor}[thm]{Corollary}
\theoremstyle{definition}

\newtheorem{exa}[thm]{Example}
\theoremstyle{remark}

\newcommand{\Fer}{\mathrm{Fer}}

\newcommand{\eqdef}{\stackrel{\rm def}{=}}

\usepackage{babel}

\newcommand{\St}{\mathcal {ST}}

\newcommand{\qc}{G(r,p,q,n)}

\newcommand{\Irr}{\mathrm{Irr}}
\newcommand{\GCD}{\mathrm{GCD}}
\newcommand{\Inv}{\mathrm{Inv}}
\newcommand{\inv}{\mathrm{inv}}
\newcommand{\Pair}{\mathrm{Pair}}
\newcommand{\pair}{\mathrm{pair}}
\newcommand{\Fix}{\mathrm{Fix}}
\newcommand{\Supp}{\mathrm{Supp}}
\author{Fabrizio Caselli}
\title[Involutory reflection groups]{Involutory reflection groups and their models}
\address{Dipartimento di matematica, Universit\`a di Bologna, \\Piazza di Porta San Donato 5, \\ Bologna 40126, Italy}

\begin{document}

\maketitle
\begin{abstract}
A finite subgroup of $GL(n,\mathbb C)$ is involutory if the sum of the dimensions of its irreducible complex representations is given by the number of absolute involutions in the group. A uniform combinatorial model is constructed for all non-exceptional irreducible complex reflection groups which are involutory including, in particular, all infinite families of finite irreducible Coxeter groups.
\end{abstract}
\section{Introduction}
In their paper \cite{BGG} Bernstein, Gelfand and Gelfand introduced the problem of the construction of a model of a group $G$, i.e. a representation which is the direct sum of all irreducible complex representations of $G$ with multiplicity one. We can find several constructions of models in the literature for the symmetric group \cite{A^+R,AA,IRS,K,K2,KV} and for some other special classes of complex reflection groups \cite{A^+R1, A,AB,B}.\\
A complex reflection group, or simply a reflection group, is a subgroup of $GL(V)$, where $V$ is a finite dimensional complex vector space, generated by reflections, i.e. by elements of finite order which fix a hyperplane pointwise. There is a well-known classification of irreducible reflection groups groups due to Chevalley \cite{C} and Shepard-Todd \cite{ST} including an infinite family $G(r,p,n)$ depending on 3 parameters together with 34 exceptional cases. 
As mentioned above one can find in the literature models some reflection groups such as the wreath product groups $G(r,1,n)$ as well as the the groups $G(2,2,n)$, which are better known as the Weyl groups of type $D$. \\
If $G$ is a finite subgroup of $GL(V)$, a theorem of Bump and Ginzburg \cite{BG} gives a combinatorial description of the character of a $G$-model if its dimension is given by the number of absolute involutions of $G$ (i.e. elements $g\in G$ such that $g\bar g=1$). We say that a group satisfying this condition is involutory. We show that a complex reflection group is involutory if and only if $\GCD(p,n)=1,2$, and the main result of this paper is an explicit and uniform construction of a model for all these groups. This construction involves in a crucial way the theory of projective reflection groups developed in \cite{Ca1}. Indeed a byproduct of this construction is also a model for some related projective reflection groups. 
The specialization of our model to the symmetric group $G(1,1,n)$ coincides with the one already appearing in \cite{A^+R,KV}, while it is apparently new in all other cases.  
The paper is organized as follows. In \S\ref{notation} we collect the notation and the preliminary results which are needed. In \S\ref{wp} we face the case of wreath products $G(r,1,n)$ separately so that the reader may understand the technical ideas developed in the general construction in a more natural way. In \S\ref{iprg} we classify all projective reflection groups of the form $G(r,p,q,n)$ (see \S\ref{notation} for the definition) which are involutory. Finally in \S\ref{models} we state and prove the main result of this work which provides a model for all involutory reflection groups. The paper ends with a conjecture about a further decomposition of the model constructed into the direct sum of two natural submodules.

\section{Notation and preliminaries}\label{notation}

In this section we collect the notations that are used in this paper as well as the preliminary results that are needed.

We let $\mathbb Z $ be the set of integer numbers and $\mathbb N$ be the set of nonnegative integer numbers. For $a,b\in \mathbb Z$, with $a\leq b$ we let $[a,b]=\{a,a+1,\ldots,b\}$ and, for $n\in \mathbb N$ we let $[n]\eqdef[1,n]$. For $r\in\mathbb N$ we let $\mathbb Z_r\eqdef \mathbb Z /r\mathbb Z$. If $r\in \mathbb N$, $r>0$, we denote by $\zeta_r$ the primitive $r$-th root of unity $\zeta_r\eqdef e^{\frac{2\pi i}{r}}$. 

The main subject of this work are the complex reflection groups \cite{Sh}, or simply reflection groups, with particular attention to their combinatorial representation theory. The most important example of a complex reflection group is the group of permutations of $[n]$, or symmetric group, that we denote by $S_n$.
We know by the work of Shephard-Todd \cite{ST} that all but a finite number of irreducible reflection groups are the groups $G(r,p,n)$ that we are going to describe. If $A$ is a matrix with complex entries we denote by $|A|$ the real matrix whose entries are the absolute values of the entries of $A$.
The \emph{wreath product} groups $G(r,n)= G(r,1,n)$ are given by all $n\times n$ matrices satysfying the following conditions: 
\begin{itemize}
\item the non-zero entries are $r$-th roots of unity;
\item there is exactly one non-zero entry in every row and every column (i.e. $|A|$ is a permutation matrix).
\end{itemize}
If $p$ divides $r$ then the reflection group $G(r,p,n)$ is the subgroup of $G(r,n)$ given by all matrices $A\in G(r,n)$ such that $\frac{\det A}{\det |A|}$ is a $\frac{r}{p}$-th root of unity. \\
Following \cite{Ca1}, a \emph{projective reflection group} is a quotient of a reflection group by a scalar subgroup. Observe that a scalar subgroup of $G(r,n)$ is necessarily a cyclic group of the form $C_q=<\zeta_qI>$ of order $q$, for some $q|r$.

It is also easy to characterize all possible scalar subgroups of the groups $G(r,p,n)$: in fact the scalar matrix $\zeta_qI$ belongs to $G(r,p,n)$ if and only if $q|r$ and $pq|rn$. In this case we let $G(r,p,q,n)\eqdef G(r,p,n)/C_q$. If $G=\qc$ then the projective reflection group $G^*\eqdef G(r,q,p,n)$, where the roles of the parameters $p$ and $q$ are interchanged, is always well-defined. We say that $G^*$ is the dual of  $G$ and we refer the reader to \cite{Ca1} for the main properties of this duality. In this paper we will see another important occurrence of the relationship between a group $G$ and its dual $G^*$.

If the non-zero entry in the $i$-th row of $g\in G(r,n)$ is $\zeta_r^{z_i}$ we let $z_i(g)\eqdef z_i\in \mathbb Z_r$ and say that $z_1(g),\ldots,z_n(g)$ are the \emph{colors} of $g$. We can also note that $g$ belongs to $G(r,p,n)$ if and only if $\sum z_i(g)\equiv 0 \mod p$.
We sometimes think of an element $g\in G(r,n)$ as a colored permutation, i.e. as the map
\begin{eqnarray*}
<\zeta_r>[n]&\rightarrow &<\zeta_r>[n]\\
\zeta_r^k i &\mapsto &\zeta_r^{k+z_i(g)}|g|(i),
\end{eqnarray*}
where $<\zeta_r>[n]$ is the set of numbers of the form $\zeta_r^{k}i$ for some $k\in \mathbb Z_r$ and $i\in [n]$, and $|g|\in S_n$ is the permutation defined by $|g|(i)=j$ if $g_{i,j}\neq 0$. We may observe that an element $g\in G(r,n)$ is uniquely determined by the permutation $|g|$ and by its colors $z_i(g)$ for all $i\in[n]$.\\
A \emph{colored cycle} $c$ (of an element $g\in G(r,n)$) is an object of the form
$$
c=\left[\begin{array}{cccc}i_1&i_2&\cdots&i_d\\ \zeta_r^{z_{i_1}}i_2&\zeta_r^{z_{i_2}}i_3&\cdots&\zeta_r^{z_{i_d}}i_1\end{array}\right],
$$
where the entries $i_1,\ldots,i_d\in [n]$ are distinct and $g(i_j)=\zeta_r^{z_{i_j}}i_{j+1}$ (with $i_{d+1}=i_1$). We define the \emph{support} of $c$ by $\Supp(c)=\{i_1,\ldots,i_d\}$, the \emph{color} of $c$ by $z(c)=\sum_j z_{i_j}$ and the \emph{length} of $c$ by $\ell(c)=d$. We say that two cycles are disjoint if their supports are. The decomposition of an element in $G(r,n)$ into the product of disjoint colored cycles is then similar to the classical one for the symmetric group.
If $g\in G(r,n)$ we let $\bar g\in G(r,n)$ be the complex conjugate of $g$. We can also observe that $\bar g$ is determined by the conditions $|\bar g|=|g|$ and $z_i(\bar g)=-z_i(g)$ for all $i\in [n]$. Since the bar operator stabilizes the cyclic subgroup $C_q=<\zeta_qI>$ it is well-defined also on the projective reflection groups $\qc$.\\ 
In \cite{Ca1} we can find a parametrization of the irreducible representations of the groups $\qc$, that we briefly recall for the reader's convenience.
Given a partition $\lambda=(\lambda_1,\ldots,\lambda_l)$ of $n$, the \emph{Ferrers diagram of shape $\lambda$} is a collection of boxes, arranged in left-justified rows, with $\lambda_i$ boxes in row $i$. 
We denote by $\Fer(r,n)$ the set of $r$-tuples $(\lambda^{(0)},\ldots,\lambda^{(r-1)})$ of Ferrers diagrams such that $\sum |\lambda^{(i)}|=n$. If $\mu \in \Fer(r,n)$ we define the \emph{color} of $\mu$ by $z(\mu)=\sum_i i|\lambda^{(i)}|$ and, if $p|r$ we let $\Fer(r,p,n)\eqdef \{\mu\in \Fer(r,n):z(\mu)\equiv0 \mod p\}$. If $\mu\in \Fer(r,n)$ we denote by $\St_{\mu}$ the set of all possible fillings of the boxes in $\mu$ with all the numbers from 1 to $n$ appearing once, in such way that rows are increasing from left to right and columns are incresing from top to bottom in every single Ferrers diagram of $\mu$. We also say that $\St_{\mu}$ is the set of \emph{standard tableaux} of shape $\mu$. Moreover we let $\St(r,n)\eqdef \cup_{\mu\in \Fer(r,n)} \St_\mu$ and we similarly define $\St(r,p,n)$.\\
If $q\in \mathbb N$ is such that $q|r$ and $pq|nr$ then the cyclic group $C_q$ acts on $\Fer(r,p,n)$ and on $\St(r,p,n)$ by a shift of $r/q$ positions of its elements (see \cite[Lemma 6.1]{Ca1}). The corresponding quotient sets are denoted by $\Fer(r,p,q,n)$ and $\St(r,p,q,n)$. 
If $T\in \St(r,p,q,n)$ we denote by $\mu(T)$ its corresponding shape in $\Fer(r,p,q,n)$ and if $\mu \in \Fer(r,p,q,n)$ we let $\St_\mu \eqdef \{T\in \St(r,p,q,n): \mu(T)=\mu\}$. 

\begin{prop}\label{dimirrep}
The irreducible complex representations of $\qc$ can be parametrized by pairs $(\mu,\rho)$, where $\mu\in \Fer(r,q,p,n)$ and $\rho\in (C_p)_{\mu}$, the stabilizer of any element in the class $\mu$ by the action of $C_p$, so that the dimension of the irreducible representation indexed by $(\mu,\rho)$ is independent of $\rho$ and it is equal to $|\St_\mu|$.
\end{prop}

In \cite[\S10]{Ca1} it is explicitly shown a generalized version of the classical Robinson-Schensted correspondence \cite[\S 7.11]{Sta}) for the symmetric groups  and of the Stanton-White \cite{SW} correspondence for the wreath products $G(r,n)$, which is valid for all projective reflection groups $\qc$. We refer to this correspondence as the \emph{projective Robinson-Schensted correspondence}. We do not describe this correspondence explicitly, but we state all the properties that we need in this paper in the following result. 

\begin{thm}\label{projRS}
There exists a map
\begin{eqnarray*}
\qc&\longrightarrow &\St(r,p,q,n)\times \St(r,p,q,n) \\
g&\longmapsto &[P(g),Q(g)],
\end{eqnarray*}
satisfying the following properties:
\begin{enumerate}
\item $P(g)$ and $Q(g)$ have the same shape in $\Fer(r,p,q,n)$ for all $g\in \qc$;
\item if $P,Q\in \St(r,p,q,n)$ have the same shape $\mu$ then
$$|\{g\in \qc:P(g)=P \textrm{ and }Q(g)=Q\}|=|(C_q)_\mu|,$$
where $(C_q)_\mu$ is the stabilizer in $C_q$ of any element in the class $\mu$;
\item if $g\mapsto [(P_0,\ldots,P_{r-1}),(Q_0,\ldots,Q_{r-1})]$
then
\begin{eqnarray*}
\bar g^{-1}&\mapsto&[(Q_0,\ldots,Q_{r-1}),(P_0,\ldots,P_{r-1})]\\
\zeta_r g &\mapsto &[(P_1,\ldots,P_{r-1},P_0),(Q_1,\ldots,Q_{r-1},Q_0)].
\end{eqnarray*}
\end{enumerate}
\end{thm}

If $G$ is a finite group we let $\Irr(G)$ be the set of irreducible complex representations of $G$. If $M$ is a complex vector space and $\rho:G\rightarrow GL(M)$ is a representation of $G$ we say that the pair $(M,\rho)$ is a $G$-\emph{model} if the character $\chi_\rho$ is the sum of the characters of all irreducible representations of $G$ over $\mathbb C$, i.e. $M$ is isomorphic as a $G$-module to the direct sum of all irreducible modules of $G$ with multiplicity one. Sometimes we simply say that $M$ is a $G$-model if we do not need to know the map $\rho$ explicitly or if it is clear from the context. It is clear that two $G$-models are always isomorphic as $G$-modules, and so we can also speak about ``the'' $G$-model.
The last result in this section is a beautiful theorem of Bump and Ginzburg, which generalizes a classical theorem of Frobenius and Schur \cite{FS}, and allows us in some cases to determine the character of the model of a finite group if we know its dimension.
 \begin{thm}[\cite{BG}, Theorem 7]\label{bugi}
Let $G$ be a finite group, $\tau\in Aut(G)$ with $\tau^2=1$ and $M$ be a $G$-model. Assume that
$$
\dim(M)=\#\{g\in G:g\tau(g)=z\},
$$ where $z$ is a central element in $G$ such that $z^2=1$.Then
$$
\chi_M(g)=\#\{u\in G:u\tau(u)=gz\}.
$$
\end{thm}

\section{The special case of wreath products}\label{wp}
In this section we let $G=G(r,n)$ and $I=I(r,n)$ be the set of \emph{absolute involutions} of $G$, i.e. elements $g$ such that $g\bar g=1$. One can check that these are exactly the symmetric matrices in $G(r,n)$. It is known \cite{A^+R} and it can be easily deduced from Proposition \ref{dimirrep} and Theorem \ref{projRS}, that the dimension of a $G$-model is equal to the cardinality of $I$. So we can consider the formal vector space
$$
M\eqdef\bigoplus_{v\in I}\mathbb C C_v
$$
having a basis indexed by the absolute involutions of $G$.
In \cite{A^+R} it is shown how one can give to $M$ the structure of a $G$-model. To describe this model we need some further notation.
If $\sigma, \tau \in S_n$ with $\tau^2=1$ we let $\inv_{\tau}(\sigma)=|\{\Inv(\sigma)\cap \Pair(\tau)|$, where 
$$\Inv(\sigma)=\{\{i,j\}:(j-i)(\sigma(j)-\sigma(i))<0\}$$
and
$$
\Pair(\tau)=\{\{i,j\}:\tau(i)=j\neq i\}.
$$
If $g,v\in G(r,n)$ with $v\bar v=1$ we let $\inv_v(g)=\inv_{|v|}(|g|)$. 
If $r$ is even we also let $$B(g, v)=\left\{\begin{array}{ll} i:|v(i)|=i,z_i(v)=2k_i(v)+1 \textrm{ with }\\k_i(v)\in[0,r/2-1]\textrm{ and }k_i(v)+z_i(g)\in[r/2,r-1].
\end{array}\right\}$$

The following result is proved in \cite{A^+R}.
\begin{thm}\label{apr}
Let $\rho:G\rightarrow GL(M)$ be defined by
$$
\rho(g)C_v=\left\{\begin{array}{ll}(-1)^{\inv_v(g)}C_{gvg^t}& \textrm{if $r$ is odd}\\(-1)^{\inv_v(g)}(-1)^{\#B(g,v)}C_{gvg^t}& \textrm{if $r$ is even},
\end{array}
\right.
$$
where $g^t$ denotes the matrix transposed of $g$.
Then $(M,\rho)$ is a $G$-model.
\end{thm}

The first target of this work is to give to $M$ another structure of a model for $G$ whose definition does not depend on the parity of $r$ and that will allow us to obtain models for other (projective) reflection groups.
If $g,g'\in G(r,n)$ we let
$$
<g,g'>=\sum_iz_i(g)z_i(g')\in \mathbb Z_r.
$$
So $<g,g'>$ is a sort of a scalar product between the color vectors of $g$ and $g'$.
\begin{thm}\label{modgrn}Let $\varrho:G\rightarrow GL(M)$ be defined by
$$
\varrho(g)C_v=\zeta_r^{<g,v>}(-1)^{\inv_v(g)}C_{|g|v|g|^{-1}}.
$$
Then $(M,\varrho)$ is a $G$-model.
\end{thm}
We observe that in this model the conjugation on the basis elements depends only on $|g|$ and so we naturally have a finer decomposition of $M$ into invariant submodules that will be partially described later in this section.
Theorem \ref{modgrn} is a particular case of the main result of this paper (Theorem \ref{main}), but we prefer to prove it separately so that the reader may understand the ideas developed in this work in a more natural way. The proof of Theorem \ref{modgrn} is splitted into several steps and, as one may have already suspected, it has a certain superposition with the proof of Theorem \ref{apr} in \cite{A^+R} so that some parts of the proof will be sketched only.\\
The first target is to prove that $\varrho$ is a representation for $G$. For this we need the following technical result.
\begin{lem}\label{zgg}
Let $g,g'\in G(r,n)$ and $\sigma\in S_n$. Then
\begin{itemize}
\item $z_i(gg')=z_i(g')+z_{|g'|(i)}(g)$;
\item $z_i(\sigma g\sigma^{-1})=z_{\sigma^{-1}(i)}(g)$.
\end{itemize}
\end{lem}
\begin{proof}
Note that $z_i(g)$ is uniquely determined by the requirement $g(i)=\zeta_r^{z_i(g)}|g|(i)$.
So we compute
$$gg'(i)=\zeta_r^{z_i(g')}g(|g'|(i))=\zeta_r^{z_i(g')}\zeta_r^{z_{|g'|(i)}(g)}|g|(|g'|(i))=\zeta_r^{z_i(g')+z_{|g'|(i)}(g)}|gg'|(i)$$
and the first part is complete.
We then apply the first part to obtain
\begin{eqnarray*}
 z_i(\sigma g\sigma^{-1})&=&z_i(\sigma^{-1})+z_{\sigma^{-1}(i)}(\sigma g)\\
&=&z_i(\sigma^{-1})+z_{\sigma^{-1}(i)}(g)+z_{|g|\sigma^{-1}(i)}(\sigma)\\
&=&z_{\sigma^{-1}(i)}(g),
\end{eqnarray*} 
where we have used the fact that $z_i(\sigma)=0$ for all $i\in [n]$ and for all $\sigma \in S_n$.
\end{proof}
For notational convenience we let $\phi_g(v)=\zeta_r^{<g,v>}(-1)^{\inv_v(g)}$, where $g\in G$ and $v\in I$.
\begin{prop}\label{hom}
For all $g,h\in G$ and $v\in I$ we have $\phi_{gh}(v)=\phi_{g}(|h|v|h|^{-1})\phi_{h}(v)$ and in particular the map $\varrho:G\rightarrow GL(M)$ defined in Theorem \ref{modgrn} is a group homomorphism.
\end{prop}
\begin{proof}To prove the claim it is enough to show that
\begin{enumerate}
\item $\inv_v(gh)\equiv \inv_v(h)+\inv_{|h|v|h|^{-1}}(g)\mod 2$;
\item $\sum z_i(gh)z_i(v)= \sum z_i(h)z_i(v)+\sum z_i(g)z_i(|h|v|h|^{-1})$.
\end{enumerate}
The first part can be proven as in \cite[Definition 6.1]{A^+R}. For the second part we apply Lemma \ref{zgg} twice in the first and in the last line of the following computation
\begin{eqnarray*}
\sum z_i(gh)z_i(v)&=&\sum z_i(h)z_i(v)+\sum z_{|h|(i)}(g)z_i(v)\\
&=& \sum z_i(h)z_i(v)+\sum z_i(g)z_{|h|^{-1}(i)}(v)\\
&=& \sum z_i(h)z_i(v)+\sum z_i(g)z_i(|h|v|h|^{-1}),
\end{eqnarray*}
and the proof is complete.
\end{proof}
Next we have to show that the character $\chi_\varrho$ of the representation $\varrho$ coincides with the character of the model of $G$.
The conditions of Theorem \ref{bugi} are satisfied for $G=G(r,n)$ with $z=1$ and $\tau$ equal to the bar operator and in particular we deduce that the value of the character of the model on an element $g$ is equal to the number of \emph{absolute square roots} of $g$, i.e. the number of elements $u\in G$ such that $u\bar u=g$. The number of absolute square roots of an element can be computed as in \cite{A^+R} and we recall it because we will need it explicitly in \S\ref{models}. For this we need the reader to go through the following observations. 

Let $c$ be a colored cycle of length $d$. If $d$ is odd then $c\bar c=c'$, where $c'$ is a cycle of length $d$ such that $z(c')=0$. If $d$ is even then $c\bar c=c_1c_2$ where $c_1$ and $c_2$ are two disjoint colored cycles of length $d/2$ such that $z(c_1)+z(c_2)=0$. On the other hand, if $c'$ is a colored cycle of length $d$, with $d$ odd, and such that $z(c')=0$, then there are exactly $r$ colored cycles $c$ of length $d$ such thar $c \bar c =c'$, and if $c_1$ and $c_2$ are two disjoint colored cycles of length $d$ such that $z(c_1)+z(c_2)=0$, then there exist exactly $rd$ cycles $c$ of length $2d$ such that $c \bar c=c_1c_2$. It is therefore natural to classify absolute square roots of a given element depending on the respective cycle structures. \\
For this reason we denote by $\Pi^{2,1}(g)$ the set of partitions of the set of disjoint cycles of $g$ into singletons and pairs of cycles having the same length. 
\begin{exa}\label{example}
If $g\in G(3,9)$ is given by $$(g(1),g(2),\ldots,g(9))=(\zeta_3 4,\zeta_3^2 2, \zeta_3^0 8, \zeta_3^2 1,\zeta_3^1 5, \zeta_3^0 7,  \zeta_3^0 6,  \zeta_3^2 9,  \zeta_3 3),$$ then the colored cycles of $g$ are 
$$c_1=\left(\begin{array}{cc} 1&4\\ \zeta_3 4&\zeta_3^2 1 \end{array}\right), \,c_2=\left(\begin{array}{c}2\\ \zeta_3^2 2\end{array}\right),\, c_3=\left(\begin{array}{ccc}3&8&9\\\zeta_3^0 8 &\zeta_3^2 9&\zeta_3 3\end{array} \right),$$
$$ c_4=\left( \begin{array}{c}5\\ \zeta_3 ^1 5\end{array} \right),\, c_5=\left( \begin{array}{cc}6&7\\ \zeta_3^0 7&  \zeta_3^0 6\end{array}\right).$$
In this case we have
$\Pi^{2,1}(g)$ contains four partitions and these are
$$\pi_1=\{\{c_1,c_5\},\{c_2,c_4\},\{c_3\}\}, \pi_2=\{\{c_1\},\{c_5\},\{c_2,c_4\},\{c_3\}\},$$
$$\pi_3=\{\{c_1,c_5\},\{c_2\},\{c_4\},\{c_3\}\},\pi_4=\{\{c_1\},\{c_5\},\{c_2\},\{c_4\},\{c_3\}\}.$$
\end{exa}

If $\pi=\{s_1,\ldots,s_h\}\in \Pi^{2,1}(g)$ we let $z(s_i)$ be the sum of the colors of the (either one or two) cycles in $s_i$, $\ell(\pi)=h$ and $\pair_j(\pi)$ be the number of pairs of cycles of length $j$ in $\pi$.\\
Let $\pi=\{s_1,\ldots,s_h\}\in \Pi^{2,1}(g)$. We say that an absolute square root $u$ of $g$ is of \emph{type} $\pi$ if the following conditions are satisfied: for every $i\in [h]$, if $s_i=\{c'\}$ is a singleton then there exists a cycle $c$ of $u$ such that $c\bar c=c'$, and if $s_i=\{c_1,c_2\}$ is a pair then there exists a cycle $c$ of $u$ such that $c \bar c=c_1c_2$. From the previous observations we have that the number of absolute square roots of type $\pi=\{s_1,\ldots,s_h\}$ is zero unless all singletons of $\pi$ have odd length and $z(s_i)=0$ for all $i\in [h]$. If these conditions are satisfied then the number of absolute square roots of type $\pi$ is $r^h\prod_jj^{\pair_j(\pi)}$. This is recorded in the following result.
\begin{prop}\label{asr}
Let $\chi$ be the character of the model of $G(r,n)$ and $g\in G(r,n)$. Then
$$
\chi(g)=\sum_{\pi}r^{\ell(\pi)}\prod_jj^{\pair_j(\pi)},
$$
where the sum is taken over all partitions $\pi\in \Pi^{2,1}(g)$ having no singletons of even length and such that $z(s)=0$ for all $s\in \pi$.
\end{prop}
In Example \ref{example} the partition $\pi_1$ is the unique element in $\Pi^{2,1}(g)$ having no singletons of even length and such that $z(s)=0$ for all $s\in \pi_1$. Therefore we have $\chi(g)=r^{\ell(\pi_1)}\prod_jj^{\pair_j(\pi_1)}=3^31^12^1=54$.

So to prove Theorem \ref{modgrn} we have to show that $\chi_{\varrho}$ agrees with the character of the model described in Proposition \ref{asr}. With this in mind we need to study the set $\Fix(g)=\{v\in I: |g|v|g|^{-1}=v\}$. Again, this can be done using ideas similar to those in \cite{A^+R}. \\ 
If $s$ is any set of cycles of $g\in G(r,n)$ we let $\Supp(s)$, the support of $s$, be the union of the supports of the cycles in $s$. If $S\subset\mathbb N$ is finite  we let $G(r,S)$ be the set of functions $g:<\zeta_r>S\rightarrow <\zeta_r>S$ such that for any $i\in S$ there exists $z_i\in \mathbb Z_r$ such that $g(\zeta_r^h i)=\zeta_r^{h+z_i}j$ for some $j\in S$, and such that the map $i\mapsto |g(i)|$ is a permutation of $S$. In particular we have that $G(r,n)=G(r,[n])$. \\
Note that the set $\Fix(g)$ depends only on $|g|$ and that there is a trivial bijection between the colored cycles of $g$ and the cycles of $|g|$ preserving supports and lengths. So, depending on the point we need to stress we will consider a partition $\pi$ of the cycles of $g$ as a partition of the cycles of $|g|$ and viceversa. \\
Let $\pi\in \Pi^{2,1}(|g|)$.
For any $s\in \pi$ we define a set of absolute involutions $S^+_{s}\subset G(r,\Supp(s))$ in the following way.
\begin{itemize}
\item If $s=\{(i_1,\ldots,i_d)\}$ is a singleton and $d$ is odd we let 
$$
S^+_{s}=\bigcup_{k\in\mathbb Z_r}\{v\in G(r,\Supp(s)):v(i_j)=\zeta_r^{k} i_j \}
$$
\item if $s=\{(i_1,\ldots,i_d)\}$ is a singleton and $d$ is even we let
$$
S^+_{s}=\bigcup_{k\in\mathbb Z_r}\Big(\{v\in G(r,\Supp(s)): v(i_j)=\zeta_r^{k}i_j \}\cup \{ v\in G(r,\Supp(s)): v(i_j)=\zeta_r^{k}i_{j+\frac{d}{2}}\}\Big)
$$
\item if $s=\{(i_1,\ldots,i_d),(j_1,\ldots,j_d)\}$ is a pair of cycles of the same length we let
$$
S^+_{s}=\bigcup_{k\in\mathbb Z_r}\bigcup_{l\in \mathbb Z_d}\{v\in G(r,\Supp(s)):v(i_h)=\zeta_r^{k} j_{h+l} \textrm{ and }v(j_h)=\zeta_r^{k} i_{h-l}\}.
$$
\end{itemize}
The reason why we use the symbol $S^+$ will be clarified in \S \ref{models}. The description of the sets $S^+_s$ is summarized in Table \ref{table1}. Inside any cell of the second column in Table \ref{table1} the parameters $k\in \mathbb Z_r$ and $l\in \mathbb Z_d$ are arbitrary but fixed. 
\begin{table}
\centering
\renewcommand{\arraystretch}{2}
\begin{tabular}{|c|c|c|}
\hline
$s$ & $S^+_s$ & $|S^+_s|$\\
\hline
\hline
$(i_1,\ldots,i_d)$, $d$ odd&\hspace{-4mm}$i_j\mapsto \zeta_r^ki_j$&$r$\\
\hline 
\raisebox{-3mm}[0pt][0pt]{$(i_1,\ldots,i_d)$, $d$ even}& $i_j\mapsto \zeta_r^ki_j $&\raisebox{-3mm}[0pt][0pt]{$2r$}\\
\cline{2-2} 
& $i_j\mapsto \zeta_r^ki_{j+\frac{d}{2}}$&\\
\hline
$(i_1,\ldots,i_d)(j_1,\ldots,j_d)$& $i_h\mapsto \zeta_r^{k} j_{h+l}$ and $j_h\mapsto \zeta_r^{k} i_{h-l}$&$dr$\\
\hline
\end{tabular}
\vspace{5mm}
\caption{}\label{table1}
\end{table}
If $\pi=\{s_1,\ldots,s_h\}$ we let $\Fix_{\pi}=\{v_1\cdots v_h:v_i\in S^+_{s_i}\}\subset I$. Then, as shown in \cite{A^+R}, we have the following description of $\Fix(g)$
$$
\Fix(g)=\bigcup_{\pi\in \Pi^{2,1}(g)}\Fix_{\pi},
$$
this being a disjoint union.

The computation of $\phi_g(v)$ can also be splitted with respect to this decomposition. 
\begin{lem}\label{deco}
Let $g\in G(r,n)$ and $\pi=\{s_1,\ldots s_h\}\in \Pi^{2,1}(g)$. Then for every $v_i\in S^+_{s_i}$ 
$$\inv_{v_1\cdots v_h}(g)=\sum_i \inv_{v_i}(g_i)$$
and
$$<g,v_1\cdots v_h>=\sum_i <g_i,v_i>,$$
where $g_i\in G(r,\Supp(s_i))$ is the restriction of $g$ to $\Supp(s_i)$.
\end{lem}
\begin{proof}
This is a straightforward verification.
\end{proof}
\begin{lem}\label{oe}
Let $g\in G(r,n)$ and $\pi=\{s_1,\ldots s_h\}\in \Pi^{2,1}(g)$. Let $g_i$ be the restriction of $g$ to $\Supp(s_i)$ and $v_i\in S^+_{s_i}$. Then $\inv_{v_i}(g_i)$ is odd if and only if $s_i=\{(i_1,\ldots,i_d)\}$ is a singleton of even length and $v_i:i_j\mapsto \zeta_r^k i_{j+\frac{d}{2}}$ for some $k\in \mathbb Z_r$.
\end{lem}
\begin{proof}
This is proved in \cite[Observation 3.7]{A^+R1}. 
\end{proof} 
We are now ready to collect these results and give a complete proof of the main result of this section.

\vspace{3mm}
\noindent\emph{Proof of Theorem \ref{modgrn}.}
By Proposition \ref{asr} we have to show that for any $g\in G$
$$
\sum_{v\in \Fix(g)}\phi_g(v)=\sum_{\pi}r^{\ell(\pi)}\prod_jj^{\pair_j(\pi)},
$$
where the last sum is taken over all partitions $\pi\in \Pi^{2,1}(g)$ having no singletons of even length and such that $z(s)=0$ for all $s\in \pi$.
We split the first sum according to the decomposition $\Fix(g)=\cup_{\pi\in\Pi^{2,1}(g)}\Fix_{\pi}$.\\
Following \cite{A^+R}, if $s=\{(i_1,\ldots,i_{d})\}\in \pi$ is a singleton of even length we define an involution $\iota:S^+_{s}\rightarrow S^+_{s}$ by 
$$
(i_j\mapsto \zeta_r^k i_j)\stackrel{\iota}{\mapsto}(i_j\mapsto \zeta_r^k i_{j+\frac{d}{2}}).
$$
Now suppose that $\pi=\{s_1,\ldots,s_h\}\in \Phi^{2,1}(g)$ has a singleton of even length, say $s_1$, and let $v_1\cdots v_h\in \Fix_\pi$ with $v_i\in S^+_{s_i}$. Then, by Lemmas \ref{deco} and \ref{oe}, we have
$$
(-1)^{\inv_{v_1\cdots v_h}(g)}=-(-1)^{\inv_{\iota(v_1)v_2\cdots v_h}(g)}
$$
whilst $<g,v_1\cdots v_h>=<g,\iota(v_1)v_2\cdots v_h>$. It follows that if $\pi$ has a part which is a singleton of even length then
$$
\sum_{v\in \Fix_{\pi}}\phi_g(v)=0.
$$
So we can restrict our attention on partitions $\pi=\{s_1,\ldots,s_h\}\in \Phi^{2,1}(g)$ having no singletons of even length. In this case  $\inv_{v}(g)$ is even for all $v\in \Fix_{\pi}$ by Lemmas  \ref{deco} and \ref{oe} and so we have 
\begin{eqnarray*}
\sum_{v\in \Fix_{\pi}}\phi_g(v)&=&\sum_{v\in \Fix_{\pi}}\zeta_r^{<g,v>}\\
&=&\prod_{i=i}^{h}\Big(\sum_{v_i\in S^+_{s_i}}\zeta_r^{<g_i,v_i>}\Big).
\end{eqnarray*}
Now, by Table \ref{table1}, any element in $S^+_{s_i}$ is the scalar multiple of an element $w_i\in S^+_{s_i}$ such that $z_j(w_i)=0$ for all $j\in \Supp (s_i)$ and, on the other hand, any scalar multiple of an element $w_i$ with this property is still in $S^+_{s_i}$. So the sum $\sum_{v_i\in S^+_{s_i}}\zeta_r^{<g_i,v_i>}$ can be splitted into sums of the form $\sum_{k=1}^{r}\zeta_r^{<g_i,\zeta_r^kw_i>}$, where $w_i\in S^+_{s_i}$ is such that $z_j(w_i)=0$ for all $j\in \Supp (s_i)$. The computation of this sum gives
$$
\sum_{k=1}^{r}\zeta_r^{<g_i,\zeta_r^kw_i>}=\sum_{k=1}^r\zeta_r^{kz(s_i)}=\sum_{k=1}^r(\zeta_r^{z(s_i)})^k.
$$
 In particular we have that if $z(s_i)$ is not zero then the sum vanishes. On the other hand, if $z(s_i)=0$ for all $i$ the previous computation shows that $\phi_g(v)=1$ for all $v\in \Fix_{\pi}$. The number of elements in $\Fix_{\pi}$ can be easily deduced from Table \ref{table1}: this is $r^{\ell(\pi)}\prod_j{j^{\pair_j(\pi)}}$ and the proof is complete.\hfill $\Box$

\vspace{3mm}

We may observe that the absolute involutions of $G(r,p,n)$ span a submodule of $M$. A natural problem is to understand which representations of $G(r,n)$ appear in this submodule. The following refinement of Theorem \ref{modgrn} establishes a model for all projective reflection groups of the form $G(r,1,p,n)$ and gives a simple solution to the previous problem. 
\begin{cor}\label{r1pn}
Let $I(r,p,n)$ be the set of absolute involutions in $G(r,p,n)$ and $$M(r,p,n)=\bigoplus_{v\in I(r,p,n)}\mathbb C C_v.$$
Then $$\chi_{M(r,p,n)}=\sum_{\mu\in \Fer(r,p,n)}\chi_\mu,$$
and in particular $(M(r,p,n),\varrho)$ is a $G(r,1,p,n)$-model.
\end{cor}
\begin{proof}
The irreducible representations of $G(r,1,p,n)$ are exactly the irreducible representations of $G(r,n)$ indexed by $\mu\in \Fer(r,p,n)$, by Proposition \ref{dimirrep}. Moreover, since $G(r,1,p,n)$ is the quotient of $G(r,n)$ by $<\zeta_r^{r/p}>$, these are also the irreducible representations of $G(r,n)$ which are fixed pointwise by the scalar element $\zeta_r^{r/p}$. Now we have that
$$
\varrho(\zeta_r^{r/p})(C_v)=\zeta_r^{<\zeta_r^{r/p},v>}(-1)^{\inv_v(\zeta_r^{r/p})}C_{|\zeta_r^{r/p}|v|\zeta_r^{r/p}|^{-1}}=\zeta_r^{\frac{r}{p} z(v)}C_v.
$$
So the basis element $C_v$ is fixed by $\zeta_r^{r/p}$ if and only if $v\in G(r,p,n)$. On the other hand the dimension of a model for $G(r,1,p,n)$ is equal to $|I(r,p,n)|$ by Proposition \ref{dimirrep} and Theorem \ref{projRS} and so there are no other independent elements fixed by $\zeta_r^{r/p}$ and the proof is complete.
\end{proof}
One may ask about some further refinements of the previous corollary. For example is it true that the representations indexed by elements $\mu\in \Fer(r,n)$ satisfying $z(\mu)=k$ are afforded by the submodule of $M$ spanned by the absolute involutions $v$ satisfying $z(v)=k$?\\
An even finer result is desirable. Let $(i_1,\ldots,i_r;j_1,\ldots,j_r)$ be a $2r$-tuple of non negative integers such that $2(i_1+\cdots+i_r)+j_1+\cdots +j_r=n$. Then the absolute involutions in $G(r,n)$ having $i_k$ 2-cycles colored with $k$ and $j_k$ 1-cycles colored with $k$ span a submodule. One may conjecture that the irreducible representations of $G(r,n)$ appearing in this submodule are exactly those corresponding to the shapes of the elements of this form by the Robinson-Schensted correspondence (Theorem \ref{projRS}).
\section{Involutory projective reflection groups}\label{iprg}

In this section we start the investigation of a model for the projective reflection groups $\qc$. The main result here is the characterization of the groups $\qc$ such that the dimension of a $\qc$-model is equal to the number of absolute involutions in $\qc$. In these groups we can directly apply Theorem \ref{bugi} the obtain a combinatorial description of the character of the model.
\begin{prop}\label{dimmod}
The dimension of the model of a projective reflection group $\qc$ is equal to the number of elements $g$ in the dual group $G(r,q,p,n)$ which correspond by means of the projective Robinson-Schensted correspondence to pairs of the form $[P,P]$, for some $P\in \St(r,q,p,n)$. 
\end{prop}
\begin{proof}
By Proposition \ref{dimirrep} we have that 
$$
\sum_{\phi\in\Irr(G(r,p,q,n))}\dim \phi=\sum_{\mu\in \Fer(r,q,p,n)}|(C_p)_{\mu}||\St_\mu|,
$$ and so the result follows from the second part of Theorem \ref{projRS}.
\end{proof}
The next target is to show that absolute involutions in $G^*$ correspond to pairs of the form $[P,P]$ under the projective Robinson-Schensted correspondence, and then to characterize those groups for which the converse holds, i.e. the groups where the fact that $v\mapsto [P,P]$ implies that $v$ is an absolute involution.

If $g\in \qc$ we say that $g$ is a \emph{symmetric} element if any (equivalently every) lift of $g$ in $G(r,n)$ is a symmetric matrix. We similarly define \emph{antisymmetric} elements in $\qc$. Observe that we can have antisymmetric elements only if $r$ is even.

\begin{lem}\label{abinsyan}
Let $g\in G(r,p,q,n)$. Then $g$ is an absolute involution, i.e. $g \bar g=1$, if and only if either $g$ is symmetric or $q$ is even and $g$ is antisymmetric.
\end{lem}
\begin{proof}
If $g$ is a symmetric element in $G(r,p,n)$ then $g \bar g=1$ and so the class of $g$ is an absolute involution in $\qc$. Similarly, if $g\in G(r,p,n)$ is antisymmetric and $q$ is even then $g \bar g=-1=\zeta_r^{\frac{r}{2}}=(\zeta_r^{\frac{r}{q}})^{\frac{q}{2}}$ and so the class of $g$ is an absolute involution in $\qc$.
Now let $g\in G(r,p,n)$ be such that the class of $g$ is an absolute involution in $\qc$. This implies that $z\eqdef z_i(g)-z_{|g|(i)}(g)$ is a multiple of $r/q$ independent of $i$. Since $|g|$ is an involution we also have that $z=-z$. This happens if $z=0$, in which case $g$ is symmetric, or $r$ is even and $z=\frac{r}{2}$. Since this is a multiple of $\frac{r}{q}$ we have that $q$ is even and that $g$ is antisymmetric.
\end{proof}

\begin{lem}\label{anrs}
Let $v\in G(r,n)$ with $r$ even. Then the following are equivalent
\begin{enumerate}
\item $v\bar v=-1$;
\item $v\mapsto [(P_0,\ldots,P_{r-1}),(P_{\frac{r}{2}},\ldots,P_{r-1},P_0,\ldots,P_{\frac{r}{2}-1})]$ for some $(P_0,\ldots,P_{r-1})\in \St(r,n)$ by the Robinson-Schensted correspondence for $G(r,n)$.
\end{enumerate}
\end{lem}
\begin{proof}
By the third part of Theorem \ref{projRS}, if
$$
v\mapsto [(P_0,\ldots,P_{r-1}),(Q_0,\ldots,Q_{r-1})]
$$
under the projective Robinson-Schensted correspondence then
$$
-v=\zeta_r^{\frac{r}{2}}v\mapsto [(P_{\frac{r}{2}},\ldots,P_{r-1},P_0,\ldots,P_{\frac{r}{2}-1}),(Q_{\frac{r}{2}},\ldots,Q_{r-1},Q_0,\ldots,Q_{\frac{r}{2}-1})].
$$
Theorem \ref{projRS} also provides us 
$$
\bar v^{-1}\mapsto[(Q_0,\ldots,Q_{r-1}),(P_0,\ldots,P_{r-1})].
$$The result follows since $v \bar v=-1$ if and only if $-v= \bar v^{-1}$ and the fact that the projective Robinson-Schensted correspondence for $G(r,n)$ is injective.
\end{proof}
We denote by $I(r,p,q,n)$ the set of absolute involutions in $\qc$.
\begin{prop}\label{geq}
Let $G=\qc$. Then $$
\sum_{\phi\in \Irr(G)}\dim \phi \geq |I(r,q,p,n)|.
$$
\end{prop}
\begin{proof}
By Proposition \ref{dimmod} it is enough to show that if $v\in I(r,q,p,n)$ is an absolute involution then $v\mapsto [P,P]$, for some $P\in\St(r,q,p,n)$. If $v$ is symmetric this is clear by the third part of Theorem \ref{projRS}. So, by Lemma \ref{abinsyan}, we can assume that $p$ is even and $v$ is antisymmetric. By Lemma \ref{anrs} any lift of $v$ in $G(r,n)$ corresponds to a pair of the form $[P,Q]=[(P_0,\ldots,P_{r-1}),(P_{\frac{r}{2}},\ldots,P_{r-1},P_0,\ldots,P_{\frac{r}{2}-1})]$. Since $p$ is even we have that $P$ and $Q$ represent the same element in $\St(r,q,p,n)$ and the proof is complete.

\end{proof}
\begin{thm}\label{invol}
 Let $G=\qc$. Then $$
\sum_{\phi\in \Irr(G)}\dim \phi=|I(r,q,p,n)|
$$ if and only if either $\GCD(p,n)=1,2$, or $\GCD(p,n)=4$ and $r\equiv p\equiv q \equiv n\equiv 4 \mod 8$.
\end{thm}
\begin{proof}
 By Proposition \ref{dimmod} and (the proof of) Proposition \ref{geq} we may deduce that $\sum_{\phi\in \Irr(G)}\dim \phi>|I(r,q,p,n)|$ if and only if there exists an element $g\in G(r,q,p,n)$ which is not an absolute involution and such that $g\mapsto[P,P]$ for some $P\in \St(r,q,p,n)$. By Lemma \ref{anrs} this is the case if and only if there exists $d|p$, $d>2$ and a shape $\mu=(\lambda^{(0)},\ldots,\lambda^{(r-1)})\in \Fer(r,q,n)$ invariant under a cyclic permutation of order $d$ (i.e. a shift of $r/d$ steps on the indices). For in this case let $P\in \St_\mu$ and $P'$ be the multitableau obtained by a cyclic permutation of the tableaux in $P$ of order $d$. Then let $v\in G(r,n)$ be such that $v\mapsto [P,P']$. We have that $P$ and $P'$ represent the same class in $\St(r,q,p,n)$ but the class of $v$ is not an absolute involution in $G(r,q,p,n)$ by Lemma \ref{anrs}.\\
The integer $d$ is necessarily also a divisor of $n$ and so we can assume that $\GCD(p,n)>2$. \\
If $\GCD(p,n)$ has an odd factor we let $d$  be any odd factor of $\GCD(p,n)$. If $\GCD(p,n)$ is a power of 2 and either $n/4$ or $r/q$ is even we let $d=4$. In all these cases we can choose $\mu=(\lambda^{(0)},\ldots,\lambda^{(r-1)})$ where
$$
\lambda^{(i)}=\left\{\begin{array}{ll}1^{n/d}& \textrm{if }i\equiv 0 \mod r/d\\ \emptyset & otherwise,
\end{array}\right.
$$
and we can easily verify that $z(\mu)=\sum i |\lambda^{(i)}|=\sum_{j=0}^{d-1}\frac{rn}{d^2}j=\frac{rn(d-1)}{2d}\equiv 0 \mod q$ and so $\mu\in \Fer(r,q,n)$, and that $\mu$ is invariant under the action of the cyclic group of order $d$.

So we are left with the case $\GCD(p,n)=4$, $n/4$ and $r/q$ odd. In this case, since $p|r$ we have $4|r$, and the condition $r/q$ odd implies $4|q$. \\
If $8|r$ let $h$ be the smallest non-negative representative of the class of $-\frac{3nr}{32}$ modulo $q/4$. Note that $0<h<r/4$. Then we let $\mu=(\lambda^{(0)},\ldots,\lambda^{(r-1)})$ where
$$
\lambda^{(i)}=\left\{\begin{array}{ll}1^{n/4-1}& \textrm{if }i\equiv 0 \mod r/4\\1& \textrm{if }i\equiv h \mod r/4\\\emptyset & otherwise,
\end{array}\right.
$$
In this case we have $z(\mu)=\sum i |\lambda^{(i)}|=4h+\frac{3nr}{8}$ and this is clearly a multiple of $q$ by construction and so $\mu\in\Fer(r,q,n)$. It is also clear that $\mu$ is invariant under the action of the cyclic group of order 4.\\
The last case to deal with is when $r\equiv 4 \mod 8$, which also forces $q\equiv 4 \mod 8$ since $r/q$ is odd. Moreover, the fact $n/4$ odd is equivalent to $n\equiv 4 \mod 8$ and so the condition $pq|nr$ together with $\GCD(n,p)=4$ implies also $p\equiv 4 \mod 8$. Under these hypothesis we have to prove that there is no shape $\mu\in \Fer(r,q,n)$ invariant under a cyclic permutation of order $4$. In fact if $\mu=(\lambda^{(0)},\ldots,\lambda^{(r-1)})$ is invariant under a cyclic permutation of order $4$ then
$$
z(\mu)=\sum_{i=0}^{r-1}i|\lambda^{(i)}|=\sum_{j=0}^3\sum_{i=0}^{r/4-1}(i+jr/4)|\lambda^{(i)}|=\frac{3nr}{8}+4\sum_{i=0}^{r/4-1}i|\lambda^{(i)}|.
$$
This can not be a multiple of $q$ since $4\not |\frac{3nr}{8}$ and the proof is complete.

\end{proof}
We conclude this section by showing that a projective reflection group $G=\qc$ and its dual group $G^*$ always have the same number of absolute involutions. This fact will be the keypoint in the description of the character of the model for the groups satisfying the conditions of Theorem \ref{invol}.
\begin{prop}\label{afm}
We always have
$$
|I(r,p,q,n)|=|I(r,q,p,n)|.
$$
\end{prop}
\begin{proof}
We will show the following stronger fact. For any involution $\sigma\in S_n$ we have
\begin{equation}\label{fhk}
|\{g\in G:g \bar g=1 \textrm{ and }|g|=\sigma\}|=|\{g\in G^*:g \bar g=1 \textrm{ and }|g|=\sigma\}|.
\end{equation}
Suppose that $\sigma$ has some fixed points. Then if $g\in G$ is such that $|g|=\sigma$ and $g\bar g=1$ then necessarily $g$ is a symmetric element and we can easily prove that
$$
|\{g\in G:g \bar g=1 \textrm{ and }|g|=\sigma\}|=\frac{r^c}{pq},
$$ 
where $c$ is the number of cycles of $\sigma$ and so Equation \eqref{fhk} clearly follows in this case.\\
Now assume that $\sigma$ has no fixed points. In this case one can show that the number of symmetric elements
$$
|\{g\in G:g\textrm{ is symmetric, }g \bar g=1 \textrm{ and }|g|=\sigma\}|=\GCD(2,p)\frac{r^c}{pq}.
$$
If $q$ is even we have to consider also antisymmetric elements. In this case we can show that

$$ \begin{array}{rl}|\{g\in G:g\textrm{ is antisymmetric,}&\hspace{-2mm}g \bar g=1 \textrm{ and }|g|=\sigma\}|=\\
&\left\{\begin{array}{ll}0& \textrm{if $p$ is even and $nr/4$ is odd  }\\ \GCD(2,p)\frac{r^c}{pq}& \textrm{otherwise}\end{array}\right.
\end{array}
$$
Summing up we have that
$$
\begin{array}{rl}|\{g\in G:g \bar g=1 &\hspace{-2mm}\textrm{and }|g|=\sigma\}|=\\
 &\hspace{-4mm}\left\{\begin{array}{ll}2\frac{r^c}{pq}& \textrm{if $p$ and $q$ are even and $nr/4$ is odd  }\\ \GCD(2,p)\GCD(2,q)\frac{r^c}{pq}& \textrm{otherwise.}\end{array}\right.
\end{array}
$$
Since these expressions are symmetric in $p$ and $q$ the proof is complete.
\end{proof}
We say that a projective reflection group $G=G(r,p,q,n)$ is \emph{involutory} if the dimension of a model of $G$ is equal to the number of absolute involutions in $G$. By Proposition \ref{afm} we have that $G(r,p,q,n)$ is involutory if and only if it satisfies the conditions in Theorem \ref{invol}. \\
If we restrict our attention on standard reflection groups we may note that a group $G(r,p,n)$ is involutory if and only if $\GCD(p,n)=1,2$. In particular all infinite families of finite irreducible Coxeter groups (these are $A_n=G(1,1,n)$, $B_n=G(2,1,n)$, $D_n=G(2,2,n)$, $I_2(r)=G(r,r,2)$) are involutory. The main goal of this work is to establish a unified construction of a model for all involutory reflection groups (and the corresponding quotients). This will be an extension of the model for the wreath products described in Theorem \ref{modgrn}, where we have to take care of the antisymmetric elements in a particular way.
\section{Models}\label{models}
From the results of the previous section we have that the dimension of the model of an involutory reflection group $G$, is equal to the number of absolute involutions of $G$ and also to the number of absolute involutions of $G^*$. In this section we show how we can give the structure of a $G$-model to the formal vector space having a basis indexed by the absolute involutions in $G^*$.

Unless otherwise stated, we let $G=G(r,p,n)$ be an involutory reflection group, i.e. such that $\GCD(p,n)=1,2$. By Theorem \ref{bugi} we have that the character $\chi$ of a $G$-model is given by
$$
\chi(g)=|\{u\in G: u\bar u=g\}|,
$$
so our first step is to compute the number of absolute square roots of a given element in $G$. We already know from \S \ref{wp} how many and which are the absolute square roots of $g$ in $G(r,n)$, so we have to understand how many of these are in $G(r,p,n)$. \\
The following is a technical result which is certainly already known but we state it for future reference and we also sketch a proof of it.
\begin{lem}\label{modeq}
Let $a_1,\ldots,a_k, b\in \mathbb Z$. Then the modular equation
$$
\sum_i a_i x_i\equiv b\mod p
$$
is solvable if and only if $d=\GCD(a_1,\ldots,a_k,p)|b$ and in this case the number of solutions modulo $p$ is $p^{k-1}d$.
\end{lem}
\begin{proof}
The condition on $d$ is clearly necessary. For the sufficiency we have to solve the equivalent equation
$$
\sum_i \frac{a_i}{d} x_i\equiv \frac{b}{d}\mod \frac{p}{d}.
$$
It is known that the fact $\GCD(a_1/d,\ldots,a_k/d,p/d)=1$ implies that the coefficients $a_1/d,\ldots,a_k/d$ represent the first line of an invertible matrix $A\in GL(k,\mathbb Z_{p/d})$. So if we consider the following change of coordinates $(x_1',\ldots,x_k')^t=A(x_1,\ldots,x_k)^t$ then the equation becomes $x_1'\equiv \frac{b}{d}\mod \frac{p}{d}$. This implies that we can choose any value for $x_2',\ldots,x_{k}'$ while the value of $x_{1}'$ is uniquely determined modulo $p/d$. The result follows.
\end{proof}
Let $u\in G(r,n)$ be such that $u \bar u=g$. Let $\pi\in \Pi^{2,1}(g)$ be the type of $u$ as defined in $\S \ref{wp}$. Recall that there is a bijection between cycles of $u$ and parts of $\pi$ preserving supports and lengths. Let $c_1,\ldots,c_h$ be the cycles of $u$ and $s_1,\ldots,s_h$ be the parts of $\pi$, where $\ell(c_i)=\ell(s_i)=\ell_i$ and $\Supp(c_i)=\Supp(s_i)$ for all $i\in [h]$. Then all elements in $G(r,n)$ obtained by multiplying (the second line of ) any cycle of $u$ by a $r$-th root of unity are still absolute square roots of $g$ of type $\pi$. So we denote by $u(x_1,\ldots,x_h)$ the absolute square root of $g$ whose cycles are $(\zeta_r^{x_1}c_1,\ldots,\zeta_r^{x_h}c_h)$. How many of these $r^h$ elements are in $G(r,p,n)$? In other words how many solutions modulo $r$ does the modular equation in $h$ variables
\begin{equation}\label{equaz}
 z(u(x_1,\ldots,x_h))=z(u)+\sum_{i=1}^h\ell_ix_i\equiv 0 \mod p
\end{equation}
have?
Now, since $\sum \ell_i=n$, if there is a cycle of odd length or $\GCD(p,n)=1$ we necessarily have $\GCD(\ell_1,\ldots,\ell_h,p)=1$ and hence, by Lemma \ref{modeq} we have exactly $p^{h-1}$ solutions modulo $p$ and so the number of solutions modulo $r$ is $p^{h-1}(r/p)^h=r^h/p$.
Now assume that all cycles have even length and that $\GCD(p,n)=2$. In this case equation \eqref{equaz} have solutions if and only if $z(u)\equiv 0 \mod 2$, by Lemma \ref{modeq}. The following result, whose proof is a straightforward verification, is the keypoint to understand when this happens.
\begin{lem}\label{parit}
Let $c=\left[\begin{array}{cccc}i_1&i_2&\cdots & i_{d}\\ \zeta_r^{z_1} i_2&\zeta_r^{z_2} i_3&\cdots&\zeta_r^{z_d}i_1\end{array}\right]$, with $d$ even, be a colored cycle. Then $c\bar c=c_1c_2$ where $$c_1=\left[\begin{array}{cccc}i_1&i_3&\cdots & i_{d-1}\\ \zeta_r^{z_1-z_2} i_3&\zeta_r^{z_3-z_4} i_5&\cdots&\zeta_r^{z_{d-1}-z_d}i_1\end{array}\right]$$ and
$$c_2=\left[\begin{array}{cccc}i_2&i_4&\cdots & i_{d}\\ \zeta_r^{z_2-z_3} i_4&\zeta_r^{z_4-z_5} i_6&\cdots&\zeta_r^{z_{d}-z_1}i_2\end{array}\right].$$
In particular, for $r$ even, $z(c)\equiv z(c_1)\equiv z(c_2) \mod 2$.
\end{lem}
By Lemma \ref{parit} the parity of $z(u)$ is equal to the number of pairs of cycles of $\pi$ having an odd color. If this number is odd Equation \eqref{equaz} has no solutions, while if it is even the number of solutions modulo $p$ is $2p^{h-1}$ and hence we have $2p^{h-1}(r/p)^h=2r^h/p$ solutions modulo $r$ of Equation \eqref{equaz}. 
We can summarize these observations in the following result.
\begin{prop}\label{asrgrpn}
Let $\GCD(p,n)=1,2$ and $g\in G(r,p,n)$. Then the number of absolute square roots of $g$ in $G(r,p,n)$ is
$$
\sum_{\phi\in \Irr(G(r,p,n))}\chi_{\phi}(g)=\frac{1}{p}\sum_{\pi}\varepsilon_{\pi}r^{\ell(\pi)}\prod_jj^{pair_j(\pi)},
$$where, if $\pi=\{s_1,\ldots,s_h\}$,
$$
\varepsilon_{\pi}=\left\{\begin{array}{ll}1,&\textrm{if }\GCD(p,n,\ell(s_1),\ldots,\ell(s_h))=1\\2&
\textrm{if $\GCD(p,n,\ell(s_1),\ldots,\ell(s_h))=2$ and the number of $s\in \pi$} \\ & \textrm{such that $s$ is a pair of cycles of odd color is even,}\\0&\textrm{otherwise,}
\end{array}\right.
$$
and the sum in the right-hand side is on all partitions $\pi\in\Pi^{2,1}(g)$ such that $z(s)=0$ for all $s\in \pi$.
\end{prop}
Once we have an algebraic-combinatorial description of  the dimension  (Theorem \ref{invol}) and of the character (Proposition \ref{asrgrpn}) of a model for $G(r,p,n)$ we have two of the main ingredients of the proof of our main result. Before state it, we need one more definition. 
If $g\in G$ and $v\in G^*$ we let 
$$
u(g,v)=z_1(\tilde v)-z_{|g|^{-1}(1)}(\tilde v)\in \mathbb Z_{r},
$$ where $\tilde v$ is any lift of $v$ in $G(r,n)$. Note that since $u(g,v)$ is the difference of two colors of $\tilde v$ it is well-defined. We denote by $I(r,p,n)^*=I(r,1,p,n)$ the set of absolute involutions in $G^*$ and we recall (Lemma \ref{abinsyan}) that these elements can be either symmetric or antisymmetric.
\begin{thm}\label{main}
Let $$M(r,p,n)^*\eqdef \bigoplus_{v\in I(r,p,n)^*}\mathbb C C_v$$ and $\varrho:G(r,p,n)\rightarrow GL(M(r,p,n)^*)$ be defined by
\begin{equation} \label{action}
\varrho(g) (C_v) \eqdef \left\{\begin{array}{ll}
\zeta_r^{<g,v>} (-1)^{\inv_v(g)}C_{|g|v|g|^{-1}} & \textrm{ if $v$ is symmetric}\\
\zeta_r^{<g,v>}\zeta_r^{u(g,v)}C_{|g|v|g|^{-1}}& \textrm{ if $v$ is antisymmetric}.\end{array}\right.
 \end{equation}
Then $(M(r,p,n)^*,\varrho)$ is a $G(r,p,n)$-model.
\end{thm}
We first prove that Equation \eqref{action} defines on $M(r,p,n)^*$ the structure of a $G(r,p,n)$-module.
\begin{lem}
The map $\varrho$ is a group homomorphism.
\end{lem}
\begin{proof}
From Proposition \ref{hom} we only have to show that $u(gh,v)=u(h,v)+u(g,|h|v|h|^{-1})$. By Lemma \ref{zgg} we have
\begin{eqnarray*}
 u(h,v)+u(g,|h|v|h|^{-1})&=& z_1(\tilde v)-z_{|h^{-1}|(1)}(\tilde v)+ z_1(|h|\tilde v|h|^{-1})-z_{|g|^{-1}(1)}(|h|\tilde v|h|^{-1})\\
&=&z_1(\tilde v)-z_{|h^{-1}|(1)}(\tilde v)+z_{|h|^{-1}(1)}(\tilde v)-z_{|h|^{-1}|g|^{-1}(1)}(\tilde v)\\
&=&z_1(\tilde v)-z_{|gh|^{-1}(1)}(\tilde v)\\
&=&u(gh,v),
\end{eqnarray*}
and the proof is complete.
\end{proof}
We observe that if $\GCD(p,n)=1$ then $\varepsilon_{\pi}=1$  for all $\pi\in \Pi^{2,1}(g)$ and for all $g\in G(r,p,n)$. In particular, by Proposition \ref{asrgrpn}, the character of the model for $G(r,p,n)$ is $\frac{1}{p}$ times the character of the model of $G(r,n)$. Moreover, in this case, the elements in $I(r,p,n)^*$ are all symmetric and so the proof of Theorem \ref{main} is a slight modification of the proof of Theorem \ref{modgrn} and is therefore left to the reader. Alternatively one can extract a proof of the case $\GCD(p,n)=1$ from the more involved situation where $\GCD(p,n)=2$. \\We can also observe that, by \cite[Theorem 4.4]{Ca1}, a reflection group $G=G(r,p,n)$ is isomorphic to its dual as an abstract group precisely if $\GCD(p,n)=1$.
 
So, unless otherwise stated, from now on we always assume that $\GCD(p,n)=2$.\\
The next target is to understand which are the elements $v\in I(r,p,n)^*$ such that $\sigma v \sigma^{-1}=v$ for a given $\sigma\in S_n$.

\begin{lem}\label{ecoac}
Let $\sigma\in S_n$ and $w\in G(r,n)$ be such that $\sigma w\sigma^{-1}=\zeta_r^{k\frac{r}{p}}w$ with $k\in [p]$. Then either $k=p$ or $k=\frac{p}{2}$.
\end{lem}
\begin{proof}
Since $z(gg')=z(g)+z(g')$ for all $g,g'\in G(r,n)$ the hypothesis forces $z(\zeta_r^{k\frac{r}{p}})=0$. This implies that $nk\frac{r}{p}\equiv 0 \mod r$, i.e. $\frac{nk}{p}$ is an integer. Therefore $k$ is a multiple of $\frac{p}{\GCD(n,p)}$ and the proof is complete.
\end{proof}

We let $\Fix(g)$ be the set of elements in $G(r,n)$ whose corresponding classes in $G(r,p,n)^*$ are absolute involutions and are fixed by the conjugation by $|g|$, i.e.
$$
\Fix(g)\eqdef\{w\in G(r,n):w\bar w\in C_p \textrm{ and }|g| w |g|^{-1}\in C_p w\}.
$$
This set will allow us to compute the character of the representation $\varrho$ since 
$$
\chi_{\varrho}(g)=\frac{1}{p}\sum_{w\in \Fix(g)}\phi_g(w),
$$
where $$
\phi_g(w)=\left\{\begin{array}{ll}
\zeta_r^{<g,w>} (-1)^{\inv_w(g)} & \textrm{ if $w$ is symmetric}\\
\zeta_r^{<g,w>}\zeta_r^{u(g,w)}& \textrm{ if $w$ is antisymmetric}.
\end{array}\right.
$$
This is simply because any element in $I(r,p,n)^*$ has $p$ lifts in $G(r,n)$, so that we are counting any element exactly $p$ times.\\
Let $g\in G$ and $w\in \Fix(g)$. Then $|g| |w| |g|^{-1}=|w|$ and in particular we can apply all the results in \S 3 of \cite{A^+R1} on the relation between the cycle structures of $g$ and $w$. From this we know that $w$ determines a partition $\pi(w)\in \Pi^{2,1}(g)$ of the set of cycles of $g$ into singletons and pairs. In this partition a cycle $c$ is a singleton of $\pi(w)$ if the restriction of $|w|$ to $\Supp(c)$ is a permutation of $\Supp(c)$ and $\{c,c'\}$ is a pair of $\pi(w)$ if the restriction of $|w|$ to $\Supp(c)$ is a bijection between $\Supp(c)$ and $\Supp(c')$ (and viceversa, since $|w|$ is an involution). From this decomposition we have
$$
\Fix(g)=\bigcup_{\pi\in \Pi^{2,1}(g)}\Fix_{\pi},
$$
this being a disjoint union, where $\Fix_{\pi}=\{w\in \Fix(g):\pi(w)=\pi\}$.
We will need a further decomposition of the sets $\Fix_{\pi}$. By Lemma \ref{abinsyan}, if $w\in G(r,n)$ is such that $w \bar w\in C_p$, then necessarily either $w\bar w=1$, i.e.  $w$ is symmetric, or $w \bar w=-1$, i.e. $w$ is antisymmetric.
Furthermore, if $w\in \Fix(g)$, by Lemma \ref{ecoac}, we have that either $|g| w |g|^{-1}=w$, i.e. $|g|$ and $w$ commute, or $|g|  w |g|^{-1}= -w$, i.e. $|g|$ and $w$ anticommute.  With this in mind we let
\begin{eqnarray*}
S^+(g)&=&\{w\in \Fix(g):w \textrm{ is symmetric and commutes with $|g|$}\}\\
S^-(g)&=&\{w\in \Fix(g):w \textrm{ is symmetric and anticommutes with $|g|$}\}\\
A^+(g)&=&\{w\in \Fix(g):w \textrm{ is antisymmetric and commutes with $|g|$}\}\\
A^-(g)&=&\{w\in \Fix(g):w \textrm{ is antisymmetric and anticommutes with $|g|$}\}\\
\end{eqnarray*}
Now let $X^{\epsilon}(g)$ be one of the previous four sets (i.e. $X=S,A$ and $\epsilon=+,-$). If $\pi\in\Pi^{2,1}(g)$ we let $X^{\epsilon}_\pi=X^{\epsilon}(g)\cap \Fix_{\pi}$.
If $\pi=\{s_1,\ldots,s_h\}\in\Pi^{2,1}(g)$ and $w\in X^{\epsilon}_\pi$ then, by construction, $w=w_1\ldots w_h$, with $w_i\in G(r,\Supp(s_i))$. Then necessarily $w_i\in X^{\epsilon}(g_i)$ for all $i$, where $g_i$ is the restriction of $g$ to $\Supp(s_i)$. 
So we have the following factorization $$X^{\epsilon}_{\{s_1,\ldots,s_h\}}=\prod_i X^\epsilon_{\{s_i\}}.$$
So the ultimate pieces that we have to evaluate are the sets of the form $X^{\epsilon}_{s}\eqdef X^{\epsilon}_{\pi}$ when $\pi=\{s\}$ has only one part (which can be either a singleton or a pair). \\
The determination of all these sets proceeds as follows. One first considers all the elements $\sigma\in S_n$ such that $|g|\sigma|g|^{-1}=\sigma$ with the given cycle structure, and then tries to put colors on $\sigma$ so that the resulting element has the requested symmetry and commuting properties.
For example let $s=\{(i_1,\ldots,i_d)\}$ with $d$ even and suppose we want to compute the set $A^-_s$. If $w\in A^-_s$ then necessarily either $|w|:i_h\mapsto i_h$ for all $h$ or $|w|:i_h\mapsto i_{h+\frac{d}{2}}$ by \cite[\S3]{A^+R1} (but the reader should better figure it out by himself). The first case is not possible since $w$ is antisymmetric so we necessarily have $|w|:i_h\mapsto i_{h+\frac{d}{2}}$. It is clear that any scalar multiple of $w$ is still in $A^+_s$ so that we can assume $z_{i_1}(w)=0$. In this case $|g|$ is the single cycle $(i_1,\ldots,i_d)$ and so from the condition $w|g|=-|g|w$ we have $w|g|(i_1)=w(i_2)=-|g|w(i_1)=-|g|(i_{1+\frac{d}{2}})=-i_{2+\frac{d}{2}}$. In particular we have $z_{i_2}(w)=r/2$. A simple recursive argument then shows that $z_{i_h}(w)=(h-1)r/2$. So we have that $w:i_h\mapsto (-1)^{h-1}i_{h+\frac{d}{2}}$. This element is antisymmetric if and only if $d/2$ is odd i.e. $d\equiv2 \mod 4$. So we deduce that, if $s=\{(i_1,\ldots,i_d)\}$, with $d$ even, then
$$
A^-_{s}=\left\{\begin{array}{ll}\cup_{k}\{w\in G(r,\Supp(s)):w(i_h)=\zeta_r^k(-1)^{h-1}i_{h+\frac{d}{2}}\}&\textrm{if }d\equiv 2 \mod 4,\\ \emptyset & \textrm{if }d\equiv 0 \mod 4. \end{array}\right.
$$
With similar reasonings one can verify the following description of these sets.
So let $\pi=\{s\}$ where $s=\{(i_1,\ldots,i_d)\}$ is a singleton. Then
\begin{itemize}
 \item if $d$ is odd 
\begin{eqnarray*}
S^+_{s}&=&\bigcup_{k\in \mathbb Z_r}\{v\in G(r,\Supp(s)):v(i_h)=\zeta_r^{k}i_h \}\\
S^-_{s}&=&\emptyset\\
A^+_{s}&=&\emptyset\\
A^-_{s}&=&\emptyset
\end{eqnarray*}

\item if $d\equiv 2 \mod 4$
\begin{eqnarray*}
S^+_{s}&=&\bigcup_{k\in \mathbb Z_r}(\{v\in G(r,\Supp(s)): v(i_h)=\zeta_r^{k}i_h\}\cup\{v\in G(r,\Supp(s)): v(i_h)=\zeta_r^{k}i_{h+\frac{d}{2}}\})\\
S^-_{s}&=&\{v\in G(r,\Supp(s)) :v(i_h)=\zeta_r^{k}(-1)^hi_h\} \\
A^+_{s}&=&\emptyset\\
A^-_{s}&=&\bigcup_{k\in \mathbb Z_r}\{v\in G(r,\Supp(s)) :v(i_h)=\zeta_r^{k} (-1)^hi_{h+\frac{d}{2}}\}
\end{eqnarray*}
\item if $d\equiv 0 \mod 4$ we have
\begin{eqnarray*}
S^+_{s}&=&\bigcup_{k\in \mathbb Z_r}(\{v\in G(r,\Supp(s)):v=\zeta_r^{k}\}\cup\{v\in G(r,\Supp(s)): v(i_h)=\zeta_r^{k}i_{h+\frac{d}{2}}\})\\
S^-_{s}&=&\bigcup_{k\in \mathbb Z_r}(\{v\in G(r,\Supp(s)):v(i_h)=\zeta_r^{k} (-1)^hi_h\}\cup\{v\in G(r,\Supp(s)):v(i_h)=\zeta_r^{k} (-1)^hi_{h+\frac{d}{2}}\}) \\
A^+_{s}&=&\emptyset\\
A^-_{s}&=&\emptyset
\end{eqnarray*}
\end{itemize}
Now suppose that $\pi=\{s\}$ has one single part which is a pair of cycles of the same length, i.e. $s=\{c_1,c_2\}$, with $c_1=(i_1,\ldots,i_d)$ and $c_2=(j_1,\ldots,j_d)$. Then
\begin{itemize}
 \item if $d$ is odd
\begin{eqnarray*}
S^+_{s}&=&\bigcup_{k\in \mathbb Z_r}\bigcup_{l\in \mathbb Z_d}\{v\in G(r,\Supp(s)):v(i_h)=\zeta_r^{k} j_{h+l} \textrm{ and }v(j_h)=\zeta_r^{k} i_{h-l}\}\\
S^-_{s}&=&\emptyset\\
A^+_{s}&=&\bigcup_{k\in \mathbb Z_r}\bigcup_{l\in \mathbb Z_d}\{v\in G(r,\Supp(s)):v(i_h)=\zeta_r^{k} j_{h+l} \textrm{ and }v(j_h)=-\zeta_r^{k} i_{h-l}\}\\
A^-_{s}&=&\emptyset
\end{eqnarray*}
\item if $d$ is even
\begin{eqnarray*}
S^+_{s}&=&\bigcup_{k\in \mathbb Z_r}\bigcup_{l\in \mathbb Z_d}\{v\in G(r,\Supp(s)):v(i_h)=\zeta_r^{k} j_{h+l} \textrm{ and }v(j_h)=\zeta_r^{k} i_{h-l}\}\\
S^-_{s}&=&\bigcup_{k\in \mathbb Z_r}\bigcup_{l\in \mathbb Z_d}\{v\in G(r,\Supp(s)):v(i_h)=\zeta_r^{k} (-1)^h j_{h+l} \textrm{ and }v(j_h)=\zeta_r^{k} (-1)^{h-l} i_{h-l}\}\\
A^+_{s}&=&\bigcup_{k\in \mathbb Z_r}\bigcup_{l\in \mathbb Z_d}\{v\in G(r,\Supp(s)):v(i_h)=\zeta_r^{k} j_{h+l} \textrm{ and }v(j_h)=-\zeta_r^{k} i_{h-l}\}\\
A^-_{s}&=&\bigcup_{k\in \mathbb Z_r}\bigcup_{l\in \mathbb Z_d}\{v\in G(r,\Supp(s)):v(i_h)=\zeta_r^{k} (-1)^h j_{h+l} \textrm{ and }v(j_h)=-\zeta_r^{k} (-1)^{h-l} i_{h-l}\}
\end{eqnarray*}
\end{itemize}
In the previous description the indices of $i_1,\ldots,i_d,j_1,\ldots,j_d$ should be considered in $\mathbb Z_d$. The description of the sets $S^+_s$, $S^-_s$, $A^+_s$ and $A^-_s$, where $s$ is either a single cycle or a pair of cycles of the same length is summarized in Table \ref{rotfloat2} for the reader's convenience. In any box of the table the parameters $k\in \mathbb Z^r$ and $l\in \mathbb Z_d$ are arbitrary but fixed.
\begin{sidewaystable}
\vspace{14cm}
\renewcommand{\arraystretch}{2}
\begin{tabular}{|c|l|l|l|l|}
\hline
$s$&$S^+_s$&$S^-_s$&$A^+_s$&$A^-_s$\\
\hline
\hline
$\{(i_1,\ldots,i_d)\}$ with $d$ odd &$i_h\mapsto\zeta_r^{k}i_h$&$\emptyset$&$\emptyset$&$\emptyset$\\
\hline
\hspace{-1.8cm}$\{(i_1,\ldots,i_d)\}$  &$i_h\mapsto\zeta_r^{k}i_h$&
\raisebox{-3mm}[0pt][0pt]{$i_h\mapsto\zeta_r^{k}(-1)^hi_h$}&\raisebox{-3mm}[0pt][0pt]{$\emptyset$}&\raisebox{-3mm}[0pt][0pt]{$i_h\mapsto\zeta_r^{k} (-1)^hi_{h+\frac{d}{2}}$}\\
\cline{2-2}
\hspace{1cm}with $d\equiv 2 \mod 4$&$i_h\mapsto\zeta_r^{k}i_{h+\frac{d}{2}}$&&&\\
\hline
\hspace{-1.8cm}$\{(i_1,\ldots,i_d)\}$& $i_h\mapsto\zeta_r^{k}i_h$&$i_h\mapsto\zeta_r^{k} (-1)^hi_h$&\raisebox{-3mm}[0pt][0pt]{$\emptyset$}&\raisebox{-3mm}[0pt][0pt]{$\emptyset$}\\
\cline{2-3}
\hspace{1cm}with $d\equiv 0 \mod 4$&$i_h\mapsto\zeta_r^{k}i_{h+\frac{d}{2}}$&$i_h\mapsto \zeta_r^{k} (-1)^hi_{h+\frac{d}{2}}$&&\\
\hline
$\{(i_1,\ldots,i_d),(j_1,\ldots,j_d)\}$, &$i_h\mapsto \zeta_r^{k} j_{h+l}$ &\raisebox{-3mm}[0pt][0pt]{$\emptyset$}&$i_h\mapsto \zeta_r^{k} j_{h+l}$ &\raisebox{-3mm}[0pt][0pt]{$\emptyset$}\\
\hspace{1cm}with $d$ odd& and $j_h\mapsto\zeta_r^{k} i_{h-l}$ & & and $j_h\mapsto-\zeta_r^{k} i_{h-l}$& \\
\hline
$\{(i_1,\ldots,i_d),(j_1,\ldots,j_d)\}$,&$i_h\mapsto\zeta_r^{k} j_{h+l}$ &$ i_h\mapsto \zeta_r^{k} (-1)^h j_{h+l}$ &$ i_h\mapsto \zeta_r^{k} j_{h+l}$ &$ i_h\mapsto \zeta_r^{k} (-1)^h j_{h+l}$\\
\hspace{1cm} with $d$ even &and $j_h\mapsto\zeta_r^{k} i_{h-l}$&and $j_h\mapsto \zeta_r^{k} (-1)^{h-l} i_{h-l}$&and $j_h\mapsto-\zeta_r^{k} i_{h-l}$& and $j_h\mapsto-\zeta_r^{k} (-1)^{h-l} i_{h-l}$\\
\hline
\end{tabular}
\vspace{5mm}
\caption[]{}\label{rotfloat2}
\end{sidewaystable}

The next result shows that elements in $\Fix(g)$ that anticommute with $|g|$ give no contribution to $\chi_{\varrho}(g)$.
\begin{lem}\label{sn-an}
Let $\pi\in \Pi^{2,1}(g)$. Then
$$
\sum_{w\in S^-_{\pi}}\phi_g(w)=-\sum_{w\in A^-_{\pi}}\phi_g(w).
$$
\end{lem}
\begin{proof}
We proceed by a case by case analysis depending on the structure of $\pi$. 
\begin{itemize}
\item $\pi$ has either a singleton of odd length or a pair of cycles of odd length. In this case the result is trivial since both $S^-_{\pi}$ and $A^-_{\pi}$ are empty (see Table \ref{rotfloat2}).
 \item $\pi=\{s_1,\ldots,s_h\}$ has a singleton of length $d$ with $d\equiv 0 \mod 4$, say $s_1=\{(i_1,\ldots,i_d)\}$. In this case $A^-_\pi$ is empty (see Table \ref{rotfloat2}) so we have to show that $\sum_{w\in S^-_{\pi}}\phi_g(w)=0$. This is similar to the proof of Theorem \ref{modgrn}. We define an involution $\iota:S^-_{s_1}\rightarrow S^-_{s_1}$ by 
$$
(i_h\mapsto \zeta_r^k (-1)^h i_h)\stackrel{\iota}{\mapsto}(i_h\mapsto \zeta_r^k (-1)^h i_{h+\frac{d}{2}}),
$$
and we extend it to an involution $\iota:S^-_\pi\rightarrow S^-_{\pi}$ by $\iota(w_1\cdots w_h)=\iota(w_1)w_2\cdots w_h$, where $w_i\in S^-_{s_i}$ for all $i$. Then, by Lemmas \ref{deco} and \ref{oe} we have $\phi_g(w)=-\phi_g(\iota(w))$ and we are done.
\item $\pi=\{s_1,\ldots,s_h\}$ has only singletons of length $d$ for some $d\equiv 2 \mod 4$ and pairs of cycles of even length. In this case we consider the bijection $\psi:S^-_{\pi}\rightarrow A^-_{\pi}$ defined as follows.
If $s=\{(i_1,\ldots,i_d)(j_1,\ldots,j_d)\}$ is a pair of cycles of odd length we let $\psi:S^-_s\rightarrow A^-_s$ be defined by
$$\begin{array}{rl}\psi: &(i_h\mapsto \zeta_r^{k} (-1)^h j_{h+l} \textrm{ and }j_h \mapsto\zeta_r^{k} (-1)^{h-l} i_{h-l})\\
&\longmapsto (i_h\mapsto \zeta_r^{k} (-1)^h j_{h+l-1} \textrm{ and }j_h\mapsto -\zeta_r^{k} (-1)^{h-l+1} i_{h-l+1})
\end{array}$$
If $s=\{(i_1,\ldots,i_d)\}$ is a singleton of length $d\equiv2 \mod 4$ we let $\psi:S^-_s\rightarrow A^-_s$ be defined by
$$\psi: (i_j \mapsto \zeta_r^{k}(-1)^ji_j)\mapsto (i_j\mapsto \zeta_r^{k} (-1)^ji_{j+\frac{d}{2}})
$$
The complete map $\psi:S^-_{\pi}\rightarrow A^-_{\pi}$ is such that $\psi(w_1\cdots w_h)=\psi(w_1)\cdots \psi(w_h)$ where $w_i\in S^-_{s_i}$. We observe that if $w\in S^-_{\pi}$ then $z_i(w)=z_i(\psi(w))$ for all $i$. Moreover $\inv(g,w)=0$ in this case by Lemma \ref{oe}. Let's compute $u(g,w)$ if $w\in A^-_{\pi}$. We have
$$
u(g,w)=z_1(w)-z_{|g|^{-1}(1)}(w).
$$ 
Let $s$ be the part of $\pi$ whose support contains 1. If $s=\{(i_1,\ldots,i_d)\}$ is a singleton we can clearly assume $1=i_1$ and $|g|^{-1}(1)=i_d$. The element $w$ is such that $w(i_j)=\zeta_r^{k} (-1)^ji_{j+\frac{d}{2}}$ by Table \ref{rotfloat2} and in particular $w(i_1)=\zeta_r^k(-1)i_{1+\frac{d}{2}}$ and $w_{i_d}=\zeta_r^k(-1)^d i_{d+\frac{d}{2}}$. In particular $z_1(w)=k+\frac{r}{2}$ and $z_{i_d}(w)=k$ since $d$ is even. It follows that $u(g,w)=\frac{r}{2}$.
If $s=\{(i_1,\ldots,i_d)(j_1,\ldots,j_d)\}$ is a pair of cycles of even length  we can clearly assume that $1=i_1$ and so $|g|^{-1}(1)=i_d$. In this case $w$ is such that $w(i_h)=\zeta_r^{k} (-1)^h j_{h+l}$ by Table \ref{rotfloat2} and in particular $z_1(w)=k+\frac{r}{2}$ and $z_{i_d}(w)=k$, since $d$ is even. So, also in this case we have $u(g,w)=\frac{r}{2}$.
Hence if $w\in S^-_{\pi}$, where $\pi$ has pairs of cycles of even length and singletons of length $\equiv 2 \mod 4$, we have
\begin{eqnarray*}
\phi_g(w)&=&\zeta_{r}^{<g,w>}(-1)^{\inv_w(g)}=\zeta_{r}^{<g,w>}\\
&=&-\zeta_r^{<g,w>}\zeta_r^{\frac{r}{2}}=-\zeta_r^{<g,\psi(w)>}\zeta_r^{u(g,\psi(w))}\\
&=&-\phi_g(\psi(w)),
\end{eqnarray*}
and the proof is complete.

\end{itemize}\end{proof}

\begin{lem}\label{sigl}
Let $\pi$ be a partition of the cycles of $g$ with a singleton of even length. Then
$$
\sum_{w\in \Fix(g)}\phi_g(w)=0.
$$
\end{lem}
\begin{proof}
Let $\pi=\{s_1,\ldots,s_h\}$ with $s_1$ a singleton of even length. Since $A^+_{\pi}=\emptyset$ by Table \ref{rotfloat2} and the contributions of $S^-_\pi$ and $A^-_\pi$ cancel each others by Lemma \ref{sn-an} we only have to show that $\sum_{w\in S^+_\pi}\phi_g(w)=0$.
This can be proved with the same argument used in Theorem \ref{modgrn} and Lemma \ref{sn-an} and the proof is therefore omitted. 
\end{proof}
We are now ready to prove the main theorem of this work.

\vspace{3mm}
\noindent \emph{Proof of Theorem \ref{main}.} We have to evaluate 
$$\sum_{w\in \Fix(g)}\phi_g(w).$$ By Lemmas \ref{sn-an} and \ref{sigl} we have
$$
\sum_{w\in \Fix(g)}\phi_g(w)=\sum_{\pi}\sum_{w\in S^+_\pi\cup A^+_\pi}\phi_g(w),
$$ where the first sum in the right-hand side is on all partitions $\pi\in \Pi^{2,1}(g)$ with no singletons of even length.\\
If $\pi=\{s_1,\ldots,s_h\}$ has a singleton of odd length then, by Table \ref{rotfloat2}, $A^+_\pi=\emptyset$. So we have
\begin{eqnarray*}
\sum_{w\in S^+_\pi\cup A^+_\pi}\phi_g(w)&=&\sum_{w\in S^+_{\pi}}\phi_g(w)\\
&=&\prod_{i=i}^{h}\left(\sum_{w_i\in S^+_{s_i}}\phi_{g_i}(w_i)\right),
\end{eqnarray*}
where $g_i$ is the restriction of $g$ to $\Supp(s_i)$.
Now, by Table \ref{rotfloat2}, any element in $S^+_{s_i}$ is the scalar multiple of an element $w_i\in S^+_{s_i}$ such that $z_j(w_i)=0$ for all $j\in \Supp (s_i)$ and, on the other hand, any scalar multiple of an element $w_i$ with this property is still in $S^+_{s_i}$. So we can apply the same argument used in the proof of Theorem \ref{modgrn} to conclude that 
$$
\sum_{w_i\in S^+_{s_i}}\phi_{g_i}(w_i)=\left\{\begin{array}{ll}0&\textrm{if }z(s_i)\neq0\\|S^+_{s_i}|& \textrm{otherwise,}
\end{array}\right.
$$ and so, if $z(s_i)=0$ for all $i\in[h]$
$$
\sum_{w\in S^+_\pi\cup A^+_\pi}\phi_g(w)=\prod_{i=1}^h|S^+_{s_i}|=r^{\ell(\pi)}\prod_j{j^{\pair_j(\pi)}}
$$by Table \ref{rotfloat2}.

Now suppose that $\pi$ has no singletons. In this case we observe that if $w\in S^+_{\pi}$ then $\inv_w(g)\equiv 0\mod 2$ by Lemma \ref{oe}, and if $w\in A^+_{\pi}$ then $u(g,w)=0$ by Table \ref{rotfloat2}. In particular we have that if $w\in S^+_{\pi} \cup A^+_{\pi}$ then $\phi_g(w)=\zeta_r^{<g,w>}$.
As in the previous paragraph we can show that 
$$\sum_{w\in S^+_\pi}\phi_g(w)=\left\{\begin{array}{ll}|S^+_\pi|&\textrm{if $z(s)=0$ for all $s\in \pi$}\\0&\textrm{otherwise}. \end{array}\right.$$
Now we consider $A^+_{s_i}$, where $s_i=\{c_1,c_2\}$. Then, by Table \ref{rotfloat2}, any element in $A^+_{s_i}$ is the scalar multiple of an element $w_i\in A^+_{s_i}$ such that $z_j(w_i)=0$ if $j\in \Supp(c_1)$ and $z_j(w_i)=r/2$ if $j\in \Supp(c_2)$. If $w_i$ has this property we can consider the previous sum restricted to all scalar multiples of $w_i$ and we obtain
\begin{eqnarray*}
\sum_{k=1}^r\zeta_r^{<g_i,\zeta_r^kw_i>}&=&\sum_{k=1}^r \zeta_r^{\sum_{j\in \Supp(c_1)}z_j(g_i)k+\sum_{j\in\Supp(c_2)}z_j(g_i)(\frac{r}{2}+k)}\\
&=& \sum_{k=1}^{r}(\zeta_r^{z(s_i})^k\zeta_r^{\frac{r}{2}z(c_2)}\\
&=&(-1)^{z(c_2)}\sum_{k=1}^{r}(\zeta_r^{z(s_i})^k.
\end{eqnarray*}
So we have that
$$\sum_{w\in A^+_{s_i}}\phi_g(w)=\left\{\begin{array}{ll}(-1)^{z(c_2)}|A^+_{s_i}|&\textrm{if }z(s_i)=0\\0&\textrm{otherwise.}
\end{array}
\right.$$
In particular, since $|A^+_{\pi}|=|S^+_{\pi}|$, by Table \ref{rotfloat2}, we have that $\sum_{w\in A^+_\pi}\phi_g(w)=\sum_{w\in S^+_\pi}\phi_g(w)$ if the number of parts of $\pi$ which are pairs of cycles of odd length is even, and $\sum_{w\in A^+_\pi}\phi_g(w)=-\sum_{w\in S^+_\pi}\phi_g(w)$ if the number of parts of $\pi$ which are pairs of cycles of odd length is odd.
The proof is complete.\hfill$\Box$

\vspace{3mm}

As in the case of the wreath products $G(r,n)$ there is a natural decomposition of $M(r,p,n)^*$ into $G(r,p,n)$-submodules. In fact, if $q|r$ and $pq|rn$ we can consider the submodule $M(r,q,p,n)\subseteq M(r,p,n)^*$ spanned by all elements $C_v$ such that $v\in I(r,q,p,n)$. The next result shows that $M(r,q,p,n)$ is the sum of all irreducible representations of $G(r,p,n)$ indexed by elements $\mu\in \Fer(r,q,p,n)$.
\begin{cor}Let $\GCD(p,n)=1,2$. Then the pair $(M(r,q,p,n),\varrho)$, where 
$$
\varrho :G(r,p,q,n)\rightarrow GL(M(r,q,p,n))
$$
is defined as in Theorem \ref{main}, is a $G(r,p,q,n)$-model.
\end{cor}
\begin{proof}
The proof is very similar to that of Corollary \ref{r1pn} and is therefore omitted. 
\end{proof}
If $\GCD(p,n)=2$, there is another natural decomposition of $M(r,p,n)^*$ into 2 $G(r,p,n)$-submodules. The submodule $Sym(r,p,n)^*$ spanned by symmetric elements and the submodule $Asym(r,p,n)^*$ spanned by antisymmetric elements. Recall from Proposition \ref{dimirrep} that an irreducible representation $\mu$ of $G(r,n)$  when restricted to $G(r,p,n)$ either remains irreducible if the stabilezer $(C_p)_{\mu}$ is trivial, or splits into two irreducible representations of $G(r,p,n)$ if $(C_p)_{\mu}$ has two elements (note that there are no other possibilities since $\GCD(p,n)=2$), and that all irreducible representations of $G(r,p,n)$ are obtained in this way. The author feels that the following can be true.
\begin{conj}
Let $\chi$ be the character of $Sym(r,p,n)^*$ and $\phi$ be an irreducible representation of $G(r,n)$. If $\phi$ does not split in $G(r,p,n)$ then $<\chi,\chi_{\phi}>=1$. If $\phi$ splits into two irreducible representations $\phi^+,\phi^-$ of $G(r,p,n)$ then
$$
<\chi,\chi_{\phi^+}>=1\Longleftrightarrow<\chi,\chi_{\phi^-}>=0.
$$
\end{conj}
In other words this conjecture says that $Sym(r,p,n)^*$ is isomorphic as a $G(r,p,n)$-module to the direct sum of all unsplit irreducible representations and of exactly one of any pair of split representations.

\emph{E-mail address: }{\tt caselli@dm.unibo.it}

\begin{thebibliography}{1}
\bibitem{A^+R1}{R. Adin, A. Postnikov and Y. Roichman}, Combinatorial Gelfand models, J. Algebra 320 (2008), 1311--1325.
\bibitem{A^+R}{R. Adin, A. Postnikov and Y. Roichman}, A Gelfand model for wreath products, Israel J. Math., in press. 
\bibitem{AA}J.L. Aguado and J.O. Araujo, A Gelfand model for the symmetric group, Communications in Algebra 29 (2001), 1841--1851.
\bibitem{A} J.O. Araujo, A Gelfand model for a Weyl group of type $B_n$, Beitr\"age Algebra Geom. 44 (2003), 359--373.
\bibitem{AB} J.O. Araujo and J.J Bigeon, A Gelfand model for a Weyl group of type $D_n$ and the branching rules $D_n \hookrightarrow B_n$, J. Algebra 294 (2005), 97--116.
\bibitem{B} R.W. Baddeley, Models and involution models for wreath products and certain Weyl groups. J. London Math. Soc. 44 (1991), 55--74.
\bibitem{BGG}I.N.Bernstein, I.M. Gelfand and S.I. Gelfand, Models of representations of compact Lie groups (Russian), Funkcional. Anal. i Prilozen. 9 (1975), 61--62. 
\bibitem{BB} {E. Bagno and R. Biagioli}, Colored-descent representations of complex reflection groups $G(r,p,n)$.  Israel J. Math.  160  (2007), 317--347
\bibitem{BG}{D. Bump and D. Ginzburg}, Generalized Frobenius-Schur numbers, J.Algebra 278 (2004), 294--313 
\bibitem{Ca1}{F. Caselli}, Projective reflection groups, preprint, {\tt arXiv:0902.0684}.
\bibitem{C}{C. Chevalley}, Invariants of finite groups generated by reflections,  Amer. J. Math.  77  (1955), 778--782
\bibitem{FS}{G. Frobenius and I. Schur}, \"Uber die reellen Darstellungen de rendlichen Gruppen, S'ber. Akad. Wiss. Berlin (1906), 186--208. 
\bibitem{IRS}N.F.J. Inglis, R.W. Richardson and J. Saxl, An explicit model for the complex representations of $S_n$, Arch. Math. (Basel) 54 (1990), 258--259.
\bibitem{K} A.A. Klyachko, Models for complex representations of groups $GL(n.q)$ and Weyl groups (Russian), Dokl. Akad. Nauk SSSR 261 (1981), 275--278.
\bibitem{K2}A.A. Klyachko, Models for complex representations of groups $GL(n.q)$ (Russian), Mat. Sb. 120 (1983), 371--386.
\bibitem{KV}V. Kodiyalam and D.-N. Verma, A natural representation model for symmetric groups, {\tt arXiv:math/0402216}
\bibitem{Sh}{G. C. Shephard}, Unitary groups generated by reflections, Canadian J. Math.  5  (1953), 364--383.
\bibitem{ST}{G. C. Shephard and J. A. Todd}, Finite unitary reflection groups,  Canadian J. Math.  6  (1954), 274--304.
\bibitem{Sta}{R. P. Stanley}, Enumerative combinatorics, vol. 2, Cambridge Studies in Advanced Mathematics 62, Cambridge University Press, Cambridge, 1999.
\bibitem{SW}{D.W. Stanton and D.E. White}, A Schensted algorithm for rim hook tableaux, J. Combin. Theory Ser. A 40 (1985), 211--247.
\end{thebibliography}
\end{document}